\documentclass{article}

\usepackage{arxiv}

\usepackage[utf8]{inputenc} 
\usepackage[T1]{fontenc}    
\usepackage{hyperref}       
\usepackage{url}            
\usepackage{booktabs}       
\usepackage{amsfonts}       
\usepackage{amsmath} 
\usepackage{amsthm} 
\usepackage{mathtools}
\usepackage{amssymb}
\usepackage{nicefrac}       
\usepackage{microtype}      
\usepackage{lipsum}
\usepackage{graphicx}
\graphicspath{ {./Img/} }

\theoremstyle{plain}

\theoremstyle{definition}

\newtheorem{remark}{Remark}
\newtheorem{example}{Example}

\title{A new approach to the reaction-diffusion systems modelling}

\author{
 Maxim Nazarov\\
  Chair of Higher Math 1\\
  Moscow Institute of Electronic Technology\\
  \texttt{Nazarov-Maximilian@yandex.ru} \\
}

\begin{document}
\maketitle
\begin{abstract}
We consider a new methodology for modelling the reaction-diffusion systems based on systems of ordinary differential equations. In contrary to the specialized numerical methods like straight line method, this new methodology is positioned as a pure alternative  at the model level for partial differential equations. In its description, the new method is largely similar to the finite volume method, but unlike the latter, it uses statistical simplifications and principles of geometric probability to describe diffusion. The main objectives of this approach are to simplify the qualitative analysis of reaction-diffusion systems and to improve the efficiency of numerical model implementation. 
The first objective is successfully addressed, as it becomes possible to use the apparatus of classical dynamical systems theory for a qualitative analysis of model dynamics based on the systems of ordinary differential equations. 
The second objective is only partially addressed, as the gain in efficiency while maintaining acceptable accuracy for numerical implementation will be significant only for certain simple initial distribution of molecules and for specific diffusion coefficients.
Furthermore, to formulate criteria for practical applicability, we separately evaluate the modelling error using this new methodology.
\end{abstract}


\section*{Introduction}

To model reaction-diffusion the default approach is to use partial differential equations of the form $\partial{u}_{i}/\partial t = D\cdot \nabla^{2} u_{i} + f(u_{1},...,u_{n})$. In most cases for numerical simulation of such PDE systems the standard\footnote{For an example of the modern use of grid methods  see \cite{Shishkin-2010}.} grid methods are implemented, based on classical  difference schemes. 
However, from a practical point of view  such an approach is far from universally rational, since for a significant portion of real chemical processes the equations will be stiff  and the computational costs will be unreasonably high. 
In practice, to simplify the modelling of reaction-diffusion systems  a transition is usually made from partial derivatives to systems of ordinary first-order differential equations. 
Most often, adaptive numerical methods are used for this transition, such as the straight line method (see \cite {Makarov-Samarski} and \cite{Gujev-Kalitkin-1992}), or alternative models are constructed using the finite volume method (for more details, see \cite{Ilyin-Gurieva, Eymard-Herbin-2000, BUSSING-MURMAN-1988}).

In the work \cite{Nazarov2011} a method of transition to ordinary differential equations was considered, which in its description is close to the finite volume method, but for the formalization of diffusion it uses statistical simplifications and principles of geometric probability. As part of the current work we will refine the description of the method from the work \cite{Nazarov2011} and estimate the error of the final model, and also consider examples of using the method to model specific reaction-diffusion systems.

\section{Description of the refined method}
\label{sec:Description}

\begin{figure}[hbt]
\centering
\includegraphics[width=0.45\textwidth]{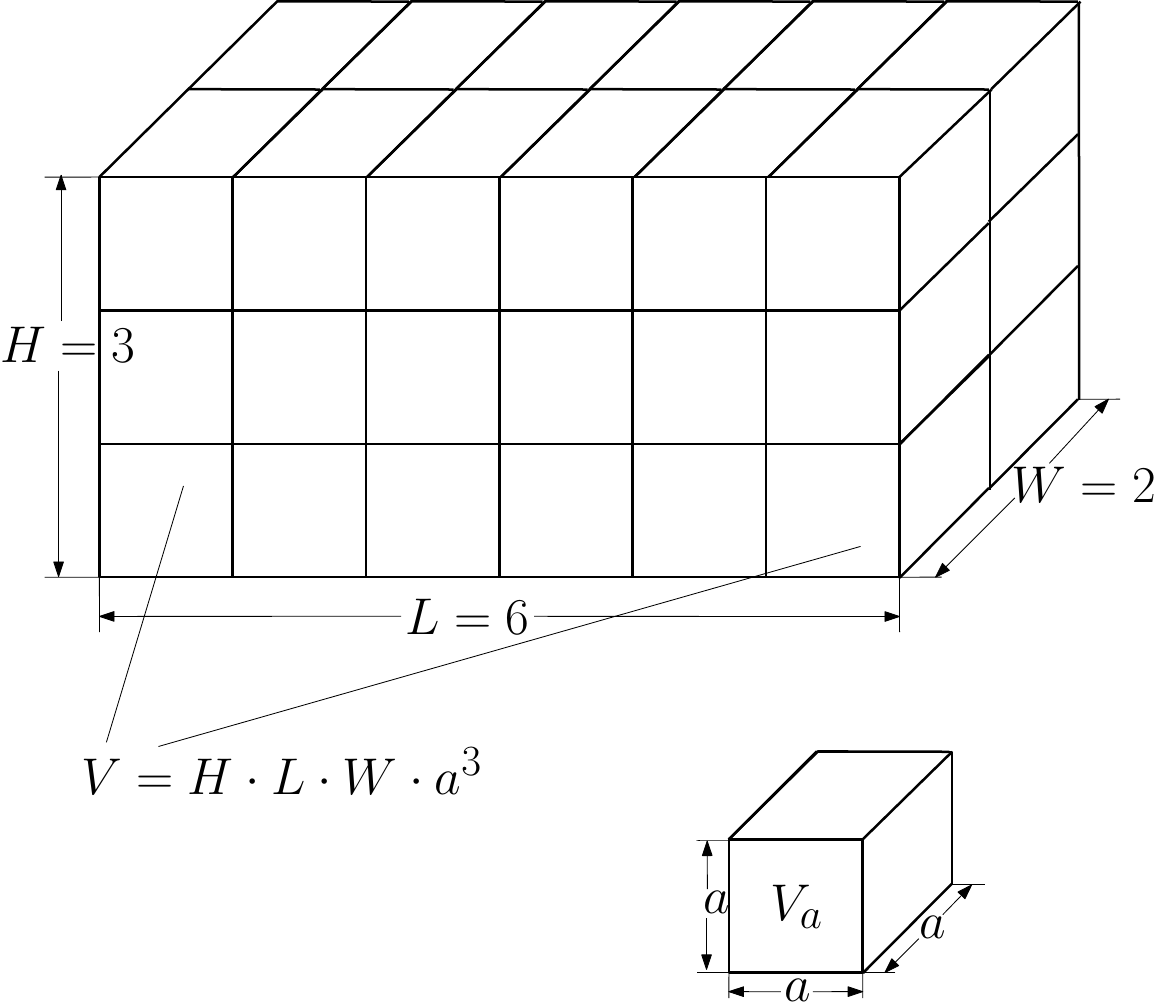}
\caption{ An example of model space partitioned into elementary cubes.}
\label{nazarov:space}
\end{figure}

As a model space  we will consider parallelepipeds of volume $ h \cdot l \cdot w $ with matched side lengths $ h = a \cdot H $, $ l = a \cdot L $  and $ w = a \cdot W $, where the parameters $ H, L $  and $ W $ are natural numbers. In this case, the model space can be uniquely partitioned into elementary cubes of volume $ V_{a}=a^3 $ (see Fig. \ref{nazarov:space}).
Let us introduce the following notations:
\begin{enumerate}
\item $ d \in \left\lbrace \uparrow \, , \downarrow \,
, \rightarrow, \leftarrow, \odot, \otimes \right\rbrace $~--- diffusion directions for elementary cubes;
\item $ \overline{d} $~--- opposite\footnote{If $d=(\downarrow)$, then the opposite direction for it will be $\overline{d} = (\uparrow)$.} direction for $ d $;
\item $ d(z,x,y) $~--- neighbouring cube with $ (z,x,y) $ in direction\footnote{For example, the cube neighbouring to $(z,x,y)$ in the direction $d=(\downarrow)$ will be $ d(z,x,y) =(z-1,x,y)$.} $ d $.
\end{enumerate} 

The description of diffusion in the model will be divided into two processes: \textbf{migration} {of substance} between the elementary cubes of space and \textbf{uniform} \textbf{intermixing} of substances in the separate cubes.

To describe these two processes  we introduce two types of variables:
\begin{itemize}
	\item The numbers $N_{i}^{(z,x,y)}(t)$ of molecules uniformly mixed in the cube $(z,x,y)$ at the moment $t$;
	\item The numbers $\mathcal{\hat{N}}_{i}^{(z,x,y)}(t)$ of molecules that entered the boundary regions of the cube at time $t$.
\end{itemize}

The dynamics of the model will be described by a system of the ordinary differential equations of the general form:

\begin{gather} \label{nazarov:Chem_Kinetics_gener:main_system_eq1}
	  \dfrac{\mathrm{d}N_{i}^{(z,x,y)}}{\mathrm{d}t} = 
		F^{(z,x,y)}_{i} + \mathrm{Dif}^{(z,x,y)}_{i (-)} + \beta_{i}	 \cdot \mathcal{\hat{N}}_{i}^{(z,x,y)},
		\quad i = \overline {1,k}
\end{gather}
\begin{gather} \label{nazarov:Chem_Kinetics_gener:main_system_eq2}
	  \dfrac{\mathrm{d}\mathcal{\hat{N}}_{i}^{(z,x,y)}}{\mathrm{d}t} = \mathcal{F}^{(z,x,y)}_{i} +
		\mathrm{Dif}^{(z,x,y)}_{i (+)} - \beta_{i}	 \cdot \mathcal{\hat{N}}_{i}^{(z,x,y)},
		\quad i = \overline {1,k}
\end{gather}

In the equations \eqref{nazarov:Chem_Kinetics_gener:main_system_eq1} and \eqref{nazarov:Chem_Kinetics_gener:main_system_eq2}, the quantity $\mathrm{Dif}^{(z,x,y)}_{i (-)}$ describes the total diffusion of molecules $i$ beyond the cube $(z,x,y)$  and the quantity $\mathrm{Dif}^{(z,x,y)}_{i (+)}$ --- the total diffusion of molecules $i$ into the boundary region of the cube $(z,x,y)$ from the   neighbouring cells.
The functions $F^{(z,x,y)}_{i}$ and $\mathcal{F}^{(z,x,y)}_{i} $ define the change in molecules $i$ in the cell $(z,x,y)$ during chemical reactions (for inert substances $F_{i}=\mathcal{F}_{i}=0 $).
The parameter $\beta_{i}$ describes the intensity of mixing of molecules $i$ within the cube $(z,x,y)$ due to internal diffusion from the boundary region of the cube.

To calculate the diffusion balance of molecules $i$ in a separate cube $(z,x,y)$  we take into account the flows of molecules $i$ in all possible diffusion directions:
\begin{gather} \label{nazarov:Diff_eq3}
	\mathrm{Dif}^{(z,x,y)}_{i(+)} = \sum_{d} 
		\xi\left( \mathrm{Emig}_{\, i, \, \overline{d}}^{d(z,x,y)} - 
			\mathrm{Emig}_{\, i,\, d}^{(z,x,y)}\right)    
\end{gather}
\begin{gather} \label{nazarov:Diff_eq3-2}
  \mathrm{Dif}^{(z,x,y)}_{i(-)} = \sum_{d} 
		\mu\left( \mathrm{Emig}_{\, i, \, \overline{d}}^{d(z,x,y)} - 
			\mathrm{Emig}_{\, i,\, d}^{(z,x,y)}\right)
\end{gather}

In the formulas \eqref{nazarov:Diff_eq3} and \eqref{nazarov:Diff_eq3-2} the functions $\xi(x) = \sigma(x)\cdot x$ and $\mu(x) = \sigma(-x)\cdot x$ are expressed either through the Heaviside function $\sigma(x)$ or through the sigmoid $\sigma^{c}(x) = {1}/({1 + e^{-x/c}})$, with $c\rightarrow 0$ and $c>0$.

The flow of molecules $\mathrm{Emig}_{i,\, d}^{(z,x,y)}(t)$ from the cube $ (z,x,y) $ in the direction $ d $ can be calculated as follows (for the derivation of the formula, see the work \cite{{Nazarov2011}}):
\begin{gather} \label{nazarov:Emig_definition_eq4}
	\mathrm{Emig}_{\, i,\, d}^{(z,x,y)} = \dfrac{3 \cdot N^{(z,x,y)}_{i} \cdot r_{i} }{16\cdot t_{r}\cdot a} +
		 \dfrac{\Delta_{(z,x,y)}^{d, i}}{ 2\cdot V_{i} \cdot M^{\, d}_{(z,x,y)} \cdot t_{r}}	
\end{gather}

In the formula \eqref{nazarov:Emig_definition_eq4}  the value $r_{i} $ specifies the maximum distance that molecules $i$ can move in time $t_{r} $  and $V_{i}$~--- is the volume occupied by an individual molecule $i$.

\begin{remark}
The probability of the diffusion must be less than one and, as a result, we obtain a constraint on the parameters: $3\cdot r_{i} \leq 16\cdot t_{r}\cdot a$. If we use the constraint equation $r_{i}= \sqrt{2 \cdot D\cdot t_{r}}$, where $D$ is the diffusion coefficient, then we can write the following inequality:
\begin{gather} \label{nazarov:ri_eq}
r_{i} \geq \dfrac{3 D}{8a}
\end{gather}
\end{remark}

The parameter $M^{\, d}_{(z,x,y)} $ in the formula \eqref{nazarov:Emig_definition_eq4} is equal to the number of different types of molecules that at time $t$ are located in one of two cubes: either in $(z,x,y)$ or in the neighbouring cube $d(z,x,y)$.
\begin{gather*} 
	M^{\, d}_{(z,x,y)} = \begin{cases} 
		 1, & \mathrm{if} \quad  N^{(z,x,y)}_{i}=N^{d(z,x,y)}_{i}=0 \, \, (\forall i); \\
		\left|  \left\lbrace i \, | \, \, N^{(z,x,y)}_{i}>0 \vee N^{d(z,x,y)}_{i}>0 \right\rbrace \right| , & \mathrm{if}	\quad \mathrm{else}.
		\end{cases}
\end{gather*}
The coefficient $ \Delta_{(z,x,y)}^{d, i} $ in the formula \eqref{nazarov:Emig_definition_eq4} gives an estimate of the relative filling of the boundary regions of the cells $(z,x,y)$ and $d(z,x,y)$ by molecules of type $i$ separately from the other types:
\begin{gather} \label{nazarov:Delta_i}
	\Delta_{(z,x,y)}^{d, i} = \dfrac{\Delta_{(z,x,y)}^{d}
	 \cdot \mathcal{\hat{N}}^{(z,x,y)}_{i}}{\sum\limits_{\mathcal{\hat{N}}_{p}>0} \mathcal{\hat{N}}^{(z,x,y)}_{p}}
\end{gather}
Let $\Delta l$ denote the maximum distance that molecules can shift $\Delta l = \max r_{i}$ in time $t_{r} $.

The coefficient $ \Delta_{(z,x,y)}^{d} $ in the formula \eqref{nazarov:Delta_i} is introduced to determine the total filling of the boundary areas of the cells $(z,x,y)$ and $d(z,x,y) $ in excess of the permissible value ($\Delta l \cdot a^{2}$): 
\begin{gather} \label{nazarov:Delta_compensation}
	\Delta_{(z,x,y)}^{d} =  \xi \left(  \max\left\lbrace L_{(z,x,y)}^{d}, 
	 L_{d(z,x,y)}^{\overline{d}} \right\rbrace - \Delta l \cdot a^{2} \right)
	 \cdot
	  \mathrm{sign} \left( L_{(z,x,y)}^{d} -  L_{d(z,x,y)}^{\overline{d}}
	 					\right)  
\end{gather}

In the equation \eqref{nazarov:Delta_compensation}  the threshold function $ \xi(x) = \sigma(x)\cdot x $ is used as $ \xi $, where $ \sigma(x) $~--- is either the Heaviside function or the sigmoid function, which is approximated to  $ \sigma^{c}(x) = {1}/({1 + e^{-x/c}})$ for $c\rightarrow 0$. The second function in the equation \eqref{nazarov:Delta_compensation} is either the sign function $ \mathrm{sign} $ or the approximate function $ \mathrm{sign}_{\tan}(x)=(2/\pi)\cdot \arctan \left( x/c\right) $ for $c\rightarrow 0$.

The value $ L_{(z,x,y)}^{d} $ from the formula \eqref{nazarov:Delta_compensation} specifies the volume that molecules in the boundary region of the cube $ (z,x,y) $ would have to occupy after an exchange with the neighbouring cube in the direction $ d $. $ L_{(z,x,y)}^{d} $ can be calculated using the formula:
\begin{gather} \label{nazarov:L_zxy_definition}
	L_{(z,x,y)}^{d} = V_{L}^{(z,x,y)} -
						 \delta V_{L}^{(z,x,y)} +
						 \delta V_{L}^{d(z,x,y)}
\end{gather} 
The value $ V_{L}^{(z,x,y)} $ from the formula \eqref{nazarov:L_zxy_definition} is the current filling of the boundary region:
\begin{gather} \label{nazarov:V_L_definition}
	V_{L}^{(z,x,y)} = \sum_{i} \mathcal{\hat{N}}^{(z,x,y)}_{i}  
	 \cdot	\dfrac{V_{i}}{k(z,x,y)} 
\end{gather} 
In this expression  $k(z,x,y)$ is the number of boundary zones of the cube $(z,x,y)$. For the interior cubes of the space  this value is constant and equals $k(z,x,y) = 6$ for $H,L,W>2$.

The value $ \delta V_{L}^{(z,x,y)}(t) $ is the potential increase in the total volume of molecules from the boundary zone, which can be divided into two parts:
\begin{gather} \label{nazarov:Delta_V_L_definition}
	\delta V_{L}^{(z,x,y)} = \delta V_{F}^{(z,x,y)}(t) + \delta V_{M}^{(z,x,y)}(t)
\end{gather} 

The first term $\delta V_{F}$  in the formula \eqref{nazarov:Delta_V_L_definition} is responsible for the increase in volume as a result of chemical reactions on the scale of the boundary zone:
\begin{gather} \label{nazarov:Delta_V_L_definition}
	\delta V_{F}^{(z,x,y)} = \sum_{i}   \mathcal{F}^{(z,x,y)}_{i} \cdot \dfrac{V_{i}}{k(z,x,y)} 
\end{gather} 

The second term $\delta V_{M}^{(z,x,y)}(t)$ in \eqref{nazarov:Delta_V_L_definition} is responsible for the increase in volume as a result of the exchange of molecules with the boundary zone of the cube $d(z,x,y)$:
\begin{gather} \label{nazarov:Delta_V_L_M_definition}
	\delta V_{M}^{(z,x,y)}(t) = \sum_{i} \dfrac{3   }{16 }\cdot N^{(z,x,y)}_{i}  \cdot    \dfrac{V_{i}\cdot r_{i}}{a } 
\end{gather} 
\begin{remark}
The proof of the correctness of the coefficient $\Delta_{(z,x,y)}^{d}$ construction is similar to the proof from \cite{Nazarov2011}.
\end{remark}
 
 To estimate the parameter $\beta_{i}$ in equations \eqref{nazarov:Chem_Kinetics_gener:main_system_eq1} and \eqref{nazarov:Chem_Kinetics_gener:main_system_eq2}  we postulate that molecules that have shifted a distance of $a/2$ within the cell will be considered intermixed. Using the constraint equation $r = \sqrt{2 D t_{s_{i}}}$  we obtain an estimate of the time $t_{s_{i}} = a^{2} / (4 D)$. For definiteness, let $(1-c)$ percent of the initial $\mathcal{\hat{N}}_{i}(0)$ dissolve during this time, where $0<c\ll1$. 
Considering the case of internal diffusion alone we can write the following expression for the $\mathcal{\hat{N}}_{i}(t)$: 
\[
	\mathcal{\hat{N}}_{i}(t) = \mathcal{\hat{N}}(0) \cdot e^{-\beta_{i} t}
\]

Substituting into this formula our data $\mathcal{\hat{N}}_{i}(t_{s_{i}}) = c \cdot \mathcal{\hat{N}}(0)$ and $t_{s_{i}} = a^{2} / (4 D)$ we obtain an estimate for the  $\beta_{i}$: 

\begin{gather} \label{nazarov:mad_coefficient}
	 \beta_{i} = \dfrac{4 D}{ a^{2}} \log \left( c^{-1} \right) =  \log \left( c^{-\frac{4 D}{ a^{2}}} \right)
\end{gather} 

The practical scope of applicability of the described method is limited, on the one hand, by the condition \eqref{nazarov:ri_eq}  and on the other hand, by the requirement that the radii $r_{i}$ be several orders of magnitude smaller than the cube length $r_{i} \leq a \cdot 10^{-2}$. For these conditions to be satisfied simultaneously  the following must be true:
\begin{gather} \label{nazarov:mad_a}
	a \geq 5 \sqrt{\dfrac{3D_{i}}{2}} \quad \forall i = \overline{1,k}
\end{gather} 

By decreasing the parameter $a$  one can achieve a more accurate description of the reaction system in space, but one must take into account the inequality \eqref{nazarov:mad_a}, as well as, the fact that the elementary cubes must be significantly larger\footnote{Otherwise  it would be impossible to consider molecules as point objects.} in volume than individual molecules: $a \gg \sqrt[3]{V_{i}} $.

 

\section{Comparative analysis of the accuracy of basic and refined models with PDE model}
\label{sec:Comparative}

We first analyse the accuracy of the generalisation scheme from \cite{Nazarov2011} by comparing it with a reference partial differential reaction-diffusion model. Then we will compare the accuracy of the refined scheme \eqref{nazarov:Chem_Kinetics_gener:main_system_eq1} -- \eqref{nazarov:mad_coefficient}  with both: reference PDE model and basic scheme from \cite{Nazarov2011}.

\begin{figure}[hbt]
\centering
\includegraphics[width=0.85\textwidth]{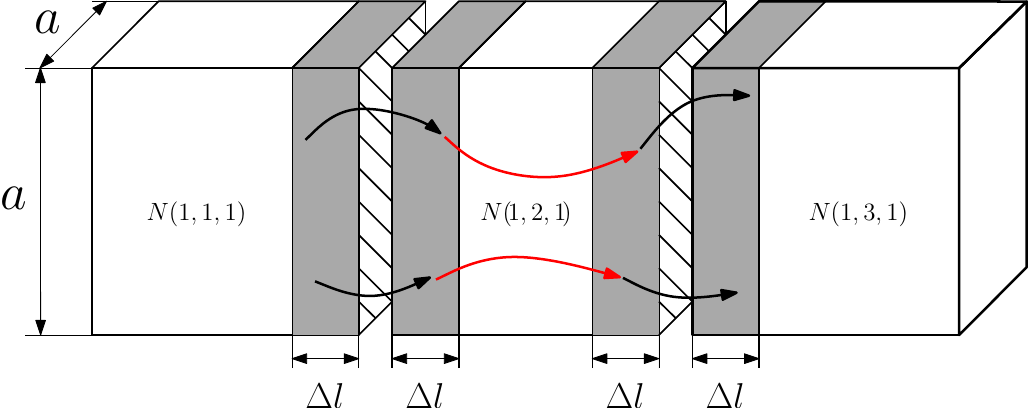}
\caption{Forwarding error: diffusion along an elementary cube in almost instant.}
\label{nazarov:forwardingerr}
\end{figure}

For simplicity  we assume that at the initial moment of time  the substance was uniformly distributed in all elementary cubes of the model space  and the substance itself was inert ($\forall i \, \, F_{i}= 0$). In this case  the main source of error for the \cite{Nazarov2011} model will be the forwarding of molecules through one cube when calculating the averaged diffusion (see figure \ref{nazarov:forwardingerr}).

The forwarding error arises because the distribution of molecules in all cubes is assumed to be uniform at all times. Thus, the model from \cite{Nazarov2011} ignores the diffusion of matter within a single cube  and any molecules entering a cube can be immediately transferred to the next cube (see the forwarding arrows in the figure \ref{nazarov:forwardingerr}).

\textbf{Problem statement 1.} To estimate the error of the model from \cite{Nazarov2011}, we will consider the simplest case:

\begin{itemize}
	\item The model has exactly one type of molecules $p$ (with parameters $V_{p}$ and $r_{p}$);
	\item At the initial moment of time molecules fill the left cube in Figure \ref{nazarov:forwardingerr}  while the other two cubes are empty.
\end{itemize}

 \[ 
	N(1,1,1) = N_{\max} = a^{3}/(V_{p}), \quad N(1,2,1) = N(1,3,1) = 0.
\]

Let us compose\footnote{We ignore all terms responsible for overflow compensation, since in the work of \cite{Nazarov2011}  when there is only one type of molecules in the model, then no overflow is possible.} system (1) from \cite{Nazarov2011} for the example in Fig. \ref{nazarov:forwardingerr}, denoting by $\alpha$ the constant $\alpha = 3\, r_{p} / (16\, t_{r}\, a )$. 
\begin{gather} \label{nazar:generSimpEq}
 \left\{\begin{aligned}
	\dfrac{\mathrm{d}N(1,1,1)}{\mathrm{d} t} &=
       	-\alpha\cdot N(1,1,1) + \alpha\cdot N(1,2,1),\\
    \dfrac{\mathrm{d}N(1,2,1)}{\mathrm{d} t} &=
       	\alpha\cdot N(1,1,1) - 2 \alpha\cdot N(1,2,1) +
       	\alpha\cdot N(1,3,1),\\
    \dfrac{\mathrm{d}N(1,3,1)}{\mathrm{d} t} &=
     	\alpha\cdot N(1,2,1) - \alpha\cdot N(1,3,1).\\
	\end{aligned}
 \right.
\end{gather}

The eigenvalues for this system of equations \ref{nazar:generSimpEq} will be $\lambda_{1}=0, \lambda_{2}=-\alpha, \lambda_{3}= -3\alpha$.
The eigenvectors that will correspond to them are $u_{1}=(1,-2,1)^{T}$, $u_{2}=(-1,0,1)^{T}$  and $u_{3}=(1,1,1)^{T}$.

The general solution of the system \ref{nazar:generSimpEq}  will take the following form (no multiple real eigenvalues):
\[
\left( \begin{array}{l}  N(1,1,1)(t)\\ N(1,2,1)(t)\\  N(1,3,1)(t) \end{array} \right)
	 = \dfrac{N_{\max}}{6}\cdot e^{-3\alpha t}\cdot \left( \begin{array}{c}  1\\  -2\\  1 \end{array} \right) -
	 \dfrac{N_{\max}}{2}\cdot e^{-\alpha t}\cdot \left( \begin{array}{c}  -1\\  0\\  1 \end{array} \right) +
	 \dfrac{N_{\max}}{3}\cdot \left( \begin{array}{c}  1\\  1\\  1 \end{array} \right)
\]

Let us now formulate the problem under consideration in terms of a partial differential model. The classical diffusion equation for the molecular density $[ N ](t)$ in the case of a constant diffusion coefficient\footnote{To estimate it  the classical equation for the average displacement radius $r(t) = \sqrt{2Dt}$ was used.} $D = r_{p}^{2}/(2 t_{r})$ takes the following form:

\[
	\dfrac{	\partial [N]}{\partial t} = 
	D \cdot \left( \dfrac{	\partial^{2} [N]}{\partial x^{2}}+ 
				\dfrac{	\partial^{2} [N]}{\partial y^{2}} +
				\dfrac{	\partial^{2} [N]}{\partial z^{2}}\right)
\]


In our case the initial conditions are $[N](z,x,y)(0) = 0$ for $x > a$ and $[N](z,x,y)(0) = N_{\max}/a^{3}$ for $x \leq a$, while  the coordinates are limited by the values $ 0 \leq z \leq a, \, 0 \leq x \leq 3a, \, 0 \leq y \leq a $. In addition, the boundary of the model space is assumed to be perfectly reflective.
The eigenfunctions of the diffusion operator for this case will be:
\[
	v_{n,m,k} = A_{n,m,k} \cdot \cos \left( \pi n \dfrac{z}{a} \right) \cdot \cos \left( \pi m \dfrac{x}{3a} \right) \cdot 
	\cos \left( \pi k \dfrac{y}{a} \right) 
\]
In turn, the eigenvalues of the diffusion operator will be the following:
\[
	\lambda_{n,m,k} = \dfrac{-D \pi^{2}}{a^{2}} \cdot 
	\left( n^{2} + m^{2}/9 + k^{2} \right)
\]
The scalar product for functions in the space under consideration is calculated as:
\[
	(f,g) = \int\limits_{0}^{a} \mathrm{d} z  \int\limits_{0}^{3a} \mathrm{d} x  
	\int\limits_{0}^{a} f(z,x,y) g(z,x,y)  \mathrm{d} y   
\]
Using the normalization condition: the scalar product $(f,g)=1$ for all eigenvectors ($v_{n,m,k}$),  one can easily find the constants $A_{n,m,k}$ (in the formula the numbers $n,m,k \geq 1$).
\[
	A_{0,0,0}\! =\! \dfrac{1}{\sqrt{3 a^3}}; \,
	A_{n,0,0}\! =\! A_{0,m,0}\! =\! A_{0,0,k}\! =\! 
	\sqrt{\dfrac{2}{3a^3}}; 
\]
\[	
	A_{n,k,0}\! =\! A_{0,m,k}\! =\! A_{n,0,m}\! 
	=\! \sqrt{\dfrac{4}{3a^3}}; \,
	A_{n,k,m}\! =\! \sqrt{\dfrac{8}{3a^3}}
\]

The general solution for the molecular density $[N]$ can be written as follows:
\[
	[N](z,x,y)(t) = \sum\limits_{n,m,k=0}^{\infty} e^{\lambda_{n,m,k}(t-t_{0})} \cdot v_{n,m,k}(z,x,y)  \int\limits_{0}^{a} \mathrm{d}z 
	\int\limits_{0}^{3a} \mathrm{d}x \int\limits_{0}^{a} 
	v_{n,m,k}(z,x,y)\cdot [N](z,x,y)(0) \mathrm{d}y
\]
Since the integrals of $\cos(\pi n\, z/a)$ and $\cos(\pi k\, y/a)$ in the range from $0$ to $a$ yield zero, only the terms for different $m$ remain in the final sum, and the remaining indices $n=k=0$. As a result, after integration, we obtain the following expression for the molecular density $[N]$:
\[
	[N](z,x,y)(t) = \dfrac{N_{\max}}{a^{3}}
	\left\lbrace \dfrac{1}{3} + \sum\limits_{m=1}^{\infty} 2 e^{\lambda_{m}t}  \cdot 
	\dfrac{\sin(\pi m/3)}{\pi m}\cdot \cos\left( \pi m \dfrac{x}{3a} \right)  \right\rbrace, \quad \lambda_{m} = -D \cdot \left( 
	\dfrac{\pi m}{3 a} \right)^{2}
\]

To estimate the forwarding error  it is necessary to obtain an exact value for the number of molecules in the third cube  in Figure \ref{nazarov:forwardingerr} (at $2a \leq x \leq 3a$). To do this we need to take the following integral:
\[
	\overline{N}(1,3,1)(t) = \int\limits_{0}^{a} \mathrm{d}z 
	\int\limits_{2a}^{3a} \mathrm{d}x \int\limits_{0}^{a} 
	 [N](z,x,y)(t) \mathrm{d}y
\]
Thus, the final expression $\overline{N}(1,3,1)(t)$ for the reference model is equal to:
\[
	\overline{N}(1,3,1)(t) = N_{\max} \left\lbrace \dfrac{1}{3} -
	 \sum\limits_{m=1}^{\infty} \dfrac{6 e^{\lambda_{m}t}}{(\pi m)^{2}}  
	 \cdot \sin(\pi m/3)\cdot \sin(2\pi m/3)
	 \right\rbrace, \quad \lambda_{m} = -D \cdot \left( 
	\dfrac{\pi m}{3 a} \right)^{2}
\]

\textbf{Evaluation of model accuracy}: The maximum value for the instantaneous error $\varepsilon_{\max}$ of the model \cite{Nazarov2011} can be estimated:
\[
	\varepsilon_{\max}  =\max_{t} \left|\overline{N}(1,3,1)(t) - {N}(1,3,1)(t)\right|
\]
Let us consider the forwarding error estimate of the model from \cite{Nazarov2011} using specific examples, choosing the values for   $N_{\max}$ and  $D$  and calculating the remaining parameters with the constraint equations.

\begin{figure}[hbt]
\centering
\includegraphics[width=0.5\textwidth]{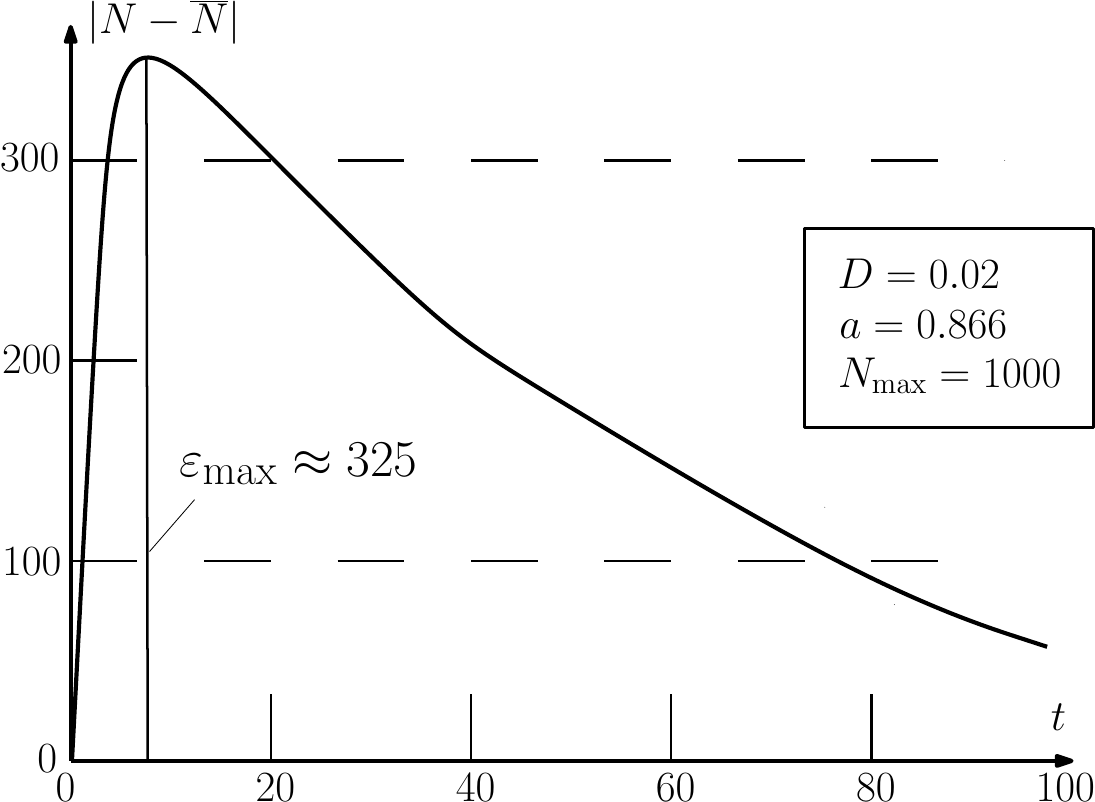}
\includegraphics[width=0.51\textwidth]{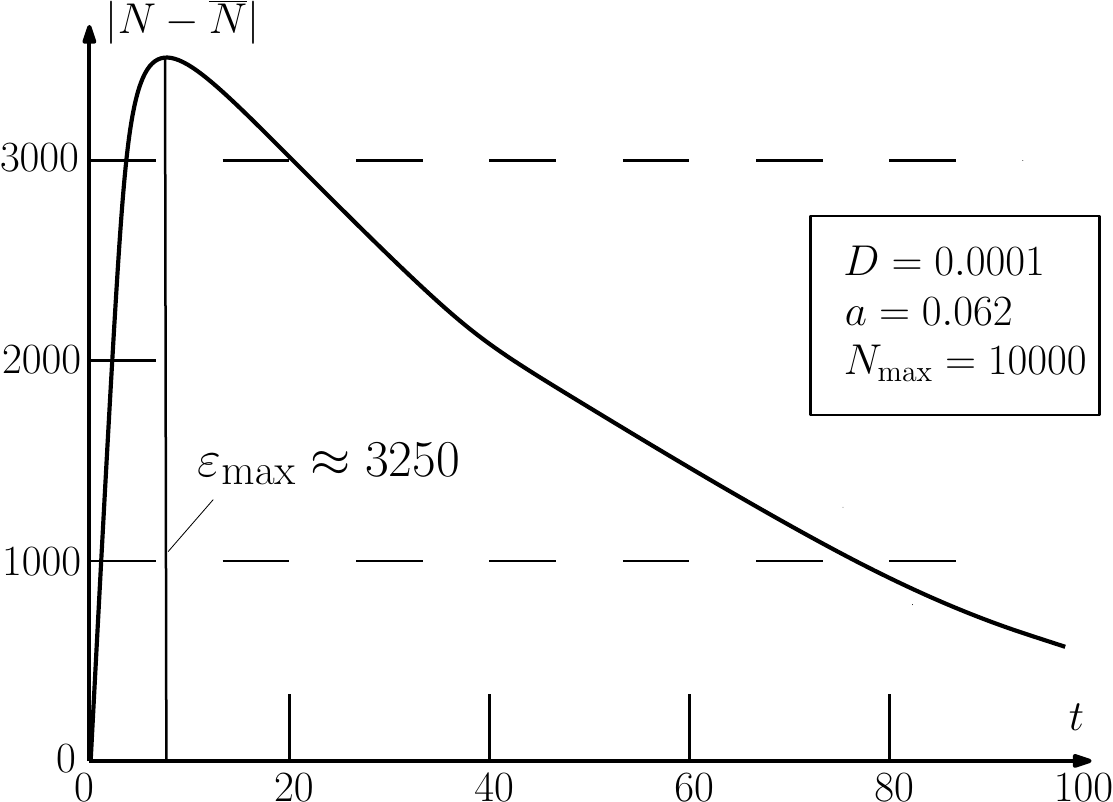}
\caption{Estimation of the mean forwarding error as a function of time for Example \ref{nazar:example1}.}
\label{nazarov:ex1}
\end{figure}

\begin{example} \label{nazar:example1}
To begin with, we set $N_{\max}=1000$ and the diffusion coefficient $D = 0.02$. If we use the constraint equations, then the length of an elementary edge will be $a \sim 5\sqrt{3D/2}$, the mean free radius $r_{p} \sim 0.01\cdot a$  and the travel time $t_{r} = r^{2}_{p}/(2D)$. As a result, we obtain the upper graph for the error $\varepsilon(t)=|\overline{N}(1,3,1)(t) - {N}(1,3,1)(t)|$ in Figure \ref{nazarov:ex1} and the value $\varepsilon_{\max}\approx 325$.
If we change the initial data to $N_{\max}=10\ 000$ and $D=0.0001$, then applying the same constraint equations for $a, r_{p}$ and $ t_{r}$ we obtain a similar graph for the error function $\varepsilon(t)$ (see the bottom half of the Figure \ref{nazarov:ex1}), and the value of the instantaneous error $\varepsilon_{\max}\approx 3250$ will increase proportionally.
\end{example}

\begin{remark}
	The relative peak error $\delta_{_{\max}} =3 \varepsilon_{_{\max}} / N_{\max}$ will be constant and independent of the choice of $N_{\max}$ and $D$, provided that the same constraint equations are used for the parameters $a, r_{p}$  and $ t_{r}$. In particular, for Example \ref{nazar:example1}   we obtain $\delta_{_{\max}} = 0.975$, which means that on short time intervals the \cite{Nazarov2011} model can be used for the case of a slightly non-uniform distribution of molecules with the constraint equations from Example  \ref{nazar:example1}. 
	However, over long time intervals with a significantly non-uniform initial distribution of molecules  the \cite{Nazarov2011} model will give an extremely poor description of diffusion, reaching a peak error of at least 97\%. 
\end{remark}

\noindent\textbf{Problem statement 2.} To estimate the error of the model \eqref{nazarov:Chem_Kinetics_gener:main_system_eq1} -- \eqref{nazarov:mad_coefficient} we postulate that:

\begin{itemize}
	\item The model involves exactly one type of molecules $p$ (with parameters $V_{p}$ and $r_{p}$);
	\item At the initial moment, molecules uniformly fill the left cell in Figure \ref{nazarov:forwardingerr}, and the other two cells are empty.
\end{itemize}
\[ 
	N(1,1,1) = N_{\max} = a^{3}/(V_{p}), \quad N(1,2,1) = N(1,3,1) = \mathcal{\hat{N}}(1,1,1) = \mathcal{\hat{N}}(1,2,1) = \mathcal{\hat{N}}(1,3,1) = 0.
\]

For this example  we cannot neglect the compensation factors $\Delta^{d}_{(z,x,y)}$, since in the \eqref{nazarov:Chem_Kinetics_gener:main_system_eq1} model -- \eqref{nazarov:mad_coefficient} they allow us to indirectly take into account the pressure in the boundary regions of the cells. 

By setting $\alpha = 3\, r_{p} / (16\, t_{r}\, a )$ and $\beta = (4 D/ a^{2} ) \log \left( c^{-1} \right)$  we can construct the following system of equations \eqref{nazarov:Chem_Kinetics_gener:main_system_eq1} and \eqref{nazarov:Chem_Kinetics_gener:main_system_eq2} for this example:
\[
 \left\{\begin{aligned}
	\dfrac{\mathrm{d}N(1,1,1)}{\mathrm{d} t} &=
       	\mu\left( \mathrm{Emig}_{\, {\leftarrow}}^{(1,2,1)} - 
			\mathrm{Emig}_{\,   \rightarrow}^{(1,1,1)}\right) + \beta\cdot \mathcal{\hat{N}}(1,1,1),\\
    \dfrac{\mathrm{d}\mathcal{\hat{N}}(1,1,1)}{\mathrm{d} t} &=
        \xi	\left( \mathrm{Emig}_{\, {\leftarrow}}^{(1,2,1)} - 
			\mathrm{Emig}_{\,   \rightarrow}^{(1,1,1)}\right) - \beta\cdot \mathcal{\hat{N}}(1,1,1),\\   	
    \dfrac{\mathrm{d}N(1,2,1)}{\mathrm{d} t} &= \mu\left( \mathrm{Emig}_{\, {\rightarrow}}^{(1,1,1)} - 
			\mathrm{Emig}_{\,   \leftarrow}^{(1,2,1)}\right) + \mu\left( \mathrm{Emig}_{\, {\leftarrow}}^{(1,3,1)} - 
			\mathrm{Emig}_{\,   \rightarrow}^{(1,2,1)}\right) 	  + \beta\cdot \mathcal{\hat{N}}(1,2,1),\\
  	\dfrac{\mathrm{d}\mathcal{\hat{N}}(1,2,1)}{\mathrm{d} t} &=
       	\xi\left( \mathrm{Emig}_{\, {\rightarrow}}^{(1,1,1)} - 
			\mathrm{Emig}_{\,   \leftarrow}^{(1,2,1)}\right) + \xi\left( \mathrm{Emig}_{\, {\leftarrow}}^{(1,3,1)} - 
			\mathrm{Emig}_{\,   \rightarrow}^{(1,2,1)}\right) - \beta\cdot \mathcal{\hat{N}}(1,2,1),\\
    \dfrac{\mathrm{d}N(1,3,1)}{\mathrm{d} t} &=
     	\mu\left( \mathrm{Emig}_{\, {\rightarrow}}^{(1,2,1)} - 
			\mathrm{Emig}_{\,  \leftarrow }^{(1,3,1)}\right) + \beta\cdot\mathcal{\hat{N}}(1,3,1),\\
     \dfrac{\mathrm{d}\mathcal{\hat{N}}(1,3,1)}{\mathrm{d} t} &=
     	\xi\left( \mathrm{Emig}_{\, {\rightarrow}}^{(1,2,1)} - 
			\mathrm{Emig}_{\,  \leftarrow }^{(1,3,1)}\right) - \beta\cdot\mathcal{\hat{N}}(1,3,1).\\ 	
	\end{aligned}
 \right.
\]

Using the formulas \eqref{nazarov:Emig_definition_eq4} -- \eqref{nazarov:mad_coefficient}  we can easily obtain expressions for all the flows  between three cubes. 
\[
	\mathrm{Emig}_{\,   \rightarrow}^{(1,1,1)} = \alpha \cdot N(1,1,1) + \dfrac{\Delta_{(1,1,1)}^{\rightarrow}}{t_{r}}; \quad 
	\mathrm{Emig}_{\,   \leftarrow}^{(1,3,1)} = \alpha \cdot N(1,3,1) + \dfrac{\Delta_{(1,3,1)}^{\leftarrow}}{t_{r}} 
\]
\[
	\mathrm{Emig}_{\,   \leftarrow}^{(1,2,1)} = \alpha \cdot N(1,2,1) + \dfrac{\Delta_{(1,2,1)}^{\leftarrow}}{t_{r}}; \quad 
	\mathrm{Emig}_{\,   \rightarrow}^{(1,2,1)} = \alpha \cdot N(1,2,1) + \dfrac{\Delta_{(1,2,1)}^{\rightarrow}}{t_{r}}
\]
\[
	\Delta_{(1,1,1)}^{\rightarrow} =  \xi \left(  \max\left\lbrace L_{(1,1,1)}^{\rightarrow}, \,	 L_{(1,2,1)}^{\leftarrow} \right\rbrace - r_{p}  \cdot a^{2} / V_{p} \right)	 \cdot	  \mathrm{sign} \left( L_{(1,1,1)}^{\rightarrow} -  L_{(1,2,1)}^{\leftarrow}	 	\right)  
\]
\[
	\Delta_{(1,2,1)}^{\rightarrow} =  \xi \left(  \max\left\lbrace L_{(1,2,1)}^{\rightarrow}, \,	 L_{(1,3,1)}^{\leftarrow} \right\rbrace - r_{p}  \cdot a^{2} / V_{p} \right)	 \cdot	  \mathrm{sign} \left( L_{(1,2,1)}^{\rightarrow} -  L_{(1,3,1)}^{\leftarrow}	 	\right)  
\]
\[
	\Delta_{(1,2,1)}^{\leftarrow} = - \Delta_{(1,1,1)}^{\rightarrow}; \quad  \Delta_{(1,3,1)}^{\leftarrow} = - \Delta_{(1,2,1)}^{\rightarrow};
\]
\[
	L_{(1,1,1)}^{\rightarrow} =  \mathcal{\hat{N}}(1,1,1) + \alpha \cdot t_{r} \cdot \left\lbrace N(1,2,1) - N(1,1,1) \right\rbrace 
\]
\[
	L_{(1,2,1)}^{\leftarrow} =   \mathcal{\hat{N}}(1,2,1)/2 + \alpha \cdot t_{r} \cdot \left\lbrace N(1,1,1) - N(1,2,1) \right\rbrace
\]
\[
	L_{(1,2,1)}^{\rightarrow} =   \mathcal{\hat{N}}(1,2,1)/2 + \alpha \cdot t_{r} \cdot \left\lbrace N(1,3,1) - N(1,2,1) \right\rbrace
\]
\[
	L_{(1,3,1)}^{\leftarrow} =   \mathcal{\hat{N}}(1,3,1)  + \alpha \cdot t_{r} \cdot \left\lbrace N(1,2,1) - N(1,3,1) \right\rbrace
\]

\textbf{Evaluation of model accuracy}: The maximum value for the instantaneous forwarding error $\varepsilon_{\max}$ of the refined model \eqref{nazarov:Chem_Kinetics_gener:main_system_eq1} -- \eqref{nazarov:mad_coefficient} can be estimated using the formula:
\[
	\varepsilon_{_{\max}} = \max_{t}\left|\overline{N}(1,3,1)(t) - N(1,3,1)(t) - \mathcal{\hat{N}}(1,3,1)(t) \right| 
\]

Let us consider the estimate $\varepsilon_{_{\max}}$ of the model \eqref{nazarov:Chem_Kinetics_gener:main_system_eq1} -- \eqref{nazarov:mad_coefficient} by choosing specific values   for the initial population $N_{\max}$ and the diffusion coefficient $D$. 

\begin{example} \label{nazar:example2}
	We use the same values  for $N_{\max}$ and $D$ as in the example \ref{nazar:example1}. Also, we use similar constraint equations for the length of the elementary edge $a \sim 5\sqrt{3D/2}$, the radius $r_{p} \sim 0.01\cdot a$  and the time $t_{r} = r^{2}_{p}/(2D)$. For the case $N_{\max}=1000$ and $D = 0.02$, we obtain the upper graph for the error $\varepsilon(t)$ in Figure \ref{nazarov:ex2}, and the value of the maximum instantaneous error will be $\varepsilon_{\max}\approx 158$. For the case $N_{\max}=10\ 000$ and $D=0.0001$ we get a similar graph below in the figure \ref{nazarov:ex2} and an error value $\varepsilon_{\max} \approx 1580$.
\begin{figure}[htb]
\centering
\includegraphics[width=0.5\textwidth]{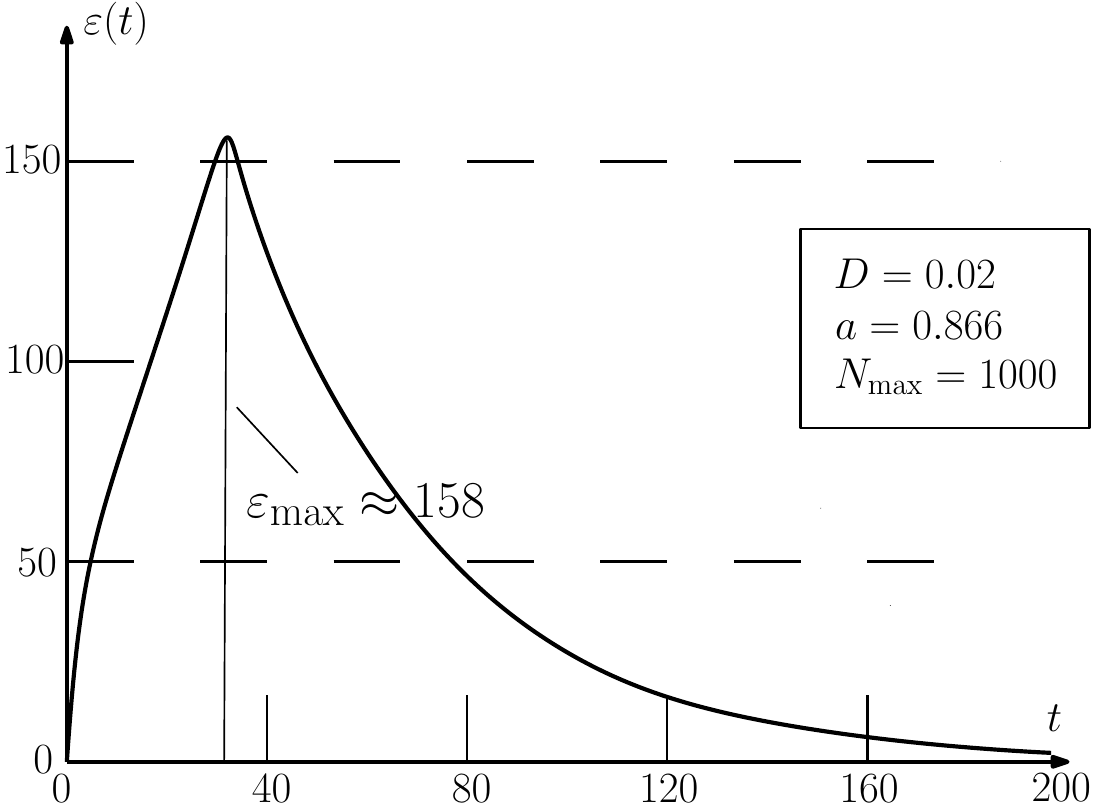}
\includegraphics[width=0.51\textwidth]{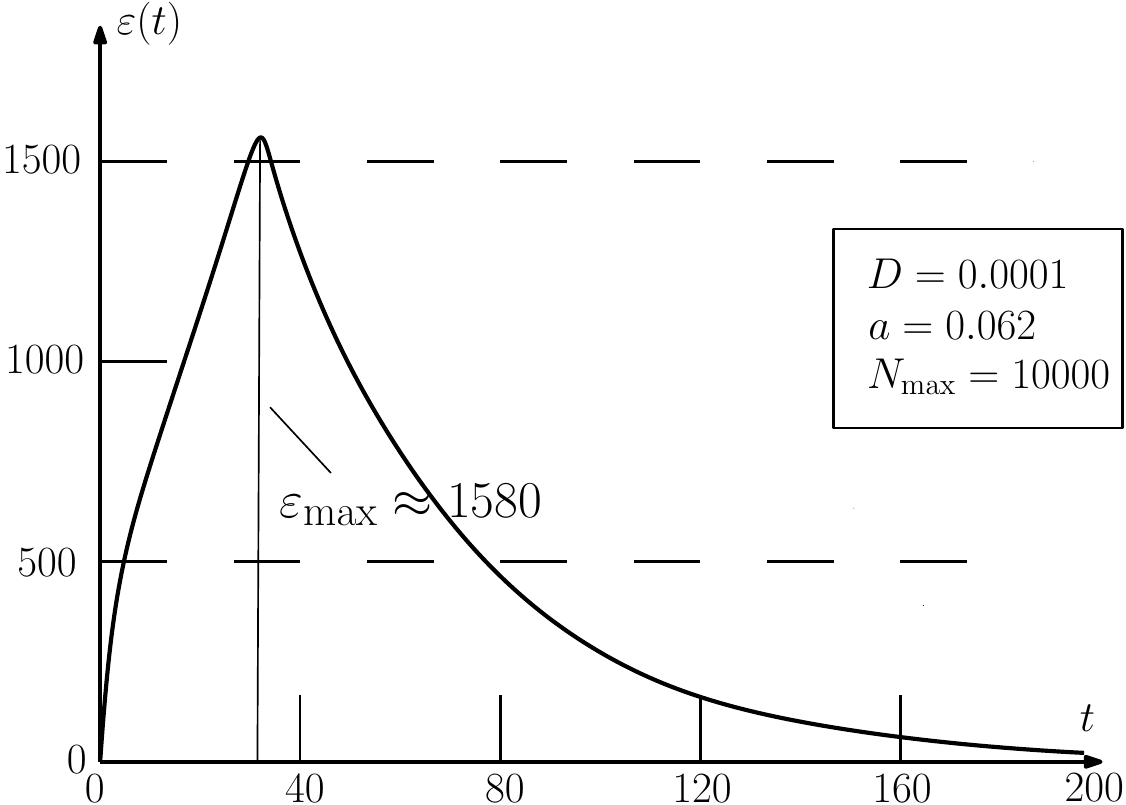}
\caption{Estimation of the mean forwarding error as a function of time for Example \ref{nazar:example2}}
\label{nazarov:ex2}
\end{figure}
\end{example}

\begin{remark}
	The relative error $\delta_{_{\max}} = \varepsilon_{_{\max}} / N_{\max}/3$ for the model \eqref{nazarov:Chem_Kinetics_gener:main_system_eq1} -- \eqref{nazarov:mad_coefficient} will also be a constant value and independent of the choice of $N_{\max}$ and $D$, provided that the same constraint equations are used for the parameters $a, r_{p}$ and $ t_{r}$. For both cases from example \ref{nazar:example2}  the value of the relative instantaneous error does not exceed $\delta_{_{\max}} = \varepsilon_{_{\max}} / N_{\max} = 0.46$, which means that the model \eqref{nazarov:Chem_Kinetics_gener:main_system_eq1} -- \eqref{nazarov:mad_coefficient} can be used\footnote{It should be noted that this fact was established exclusively for the constraint equations considered in examples \ref{nazar:example1} and \ref{nazar:example2}: $a \sim 5\sqrt{3D/2}$, $r_{p} \sim 0.01\cdot a$ and $t_{r} = r^{2}_{p}/(2D)$.} for the case of a more heterogeneous distribution of molecules over relatively longer time intervals. In fact, for this example, the \eqref{nazarov:Chem_Kinetics_gener:main_system_eq1} model -- \eqref{nazarov:mad_coefficient} achieves a peak error of 46\%, which is  less than the error of the model from \cite{Nazarov2011} of 97\%.
\end{remark}




\section*{Conclusion}

\noindent The main advantages of the \eqref{nazarov:Chem_Kinetics_gener:main_system_eq1} -- \eqref{nazarov:mad_coefficient} model compared to other approaches for describing reaction-diffusion systems are, respectively:

\begin{itemize}
	\item Possibility of analysing reaction-diffusion dynamics by using methods from dynamic systems theory.
	\item Significant gain in computational costs compared to the implementation of the model based on PDE.
\end{itemize}

\noindent The disadvantages of the reaction-diffusion model \eqref{nazarov:Chem_Kinetics_gener:main_system_eq1} -- \eqref{nazarov:mad_coefficient} are respectively:

\begin{itemize}
	\item The presence of a lower limit for the size of the elementary cube $a$, which makes the model inapplicable for the case of sufficiently small cubes;
	\item Poor integration with models that include the active flow of matter in space;
	\item Significantly lower modelling accuracy compared to partial derivative models.
\end{itemize}

One can use analogous scheme to model \eqref{nazarov:Chem_Kinetics_gener:main_system_eq1} -- \eqref{nazarov:mad_coefficient} in order to describe heat transfer instead of diffusion.  An example of such a scheme for a simpler version of the model is presented in the paper \cite{Nazarov_Chem2}. However, it should be noted that the chemical kinetics model used in \cite{Nazarov_Chem2} was not sufficiently rigorous, and therefore, for a more accurate description of reaction dynamics, it is recommended to use the chemical kinetics model from the paper \cite{Nazarov_Chem3}.

Promising generalisations for the reaction-diffusion model \eqref{nazarov:Chem_Kinetics_gener:main_system_eq1} -- \eqref{nazarov:mad_coefficient} are the search for an effective way of incorporating laminar flow of matter in model formalism, as well as, the assessment of the model error from the overlap of the boundary regions of the elementary cubes.

In conclusion, we note that the model \eqref{nazarov:Chem_Kinetics_gener:main_system_eq1} -- \eqref{nazarov:mad_coefficient} can be used not only for the modelling of chemical kinetics, but also for describing population dynamics, as well as, for modelling multicellular organisms in Biology.


\bibliographystyle{unsrt}  


\end{document}